\documentclass[11pt]{article}

\usepackage[letterpaper,left=1.25in,right=1.25in,top=1in,bottom=1in]{geometry}
\usepackage{amsmath,amssymb,amsthm}
\usepackage{setspace}
\usepackage{microtype}
\usepackage{needspace}
\usepackage{hyperref}

\allowdisplaybreaks
\newtheorem{theorem}{Theorem}[section]

\newtheorem{lemma}[theorem]{Lemma}
\newtheorem{corollary}[theorem]{Corollary}
\newtheorem{claim}{Claim}[section]
\newcommand{\del}{\backslash}
\newcommand{\tw}{\operatorname{tw}}
\newcommand{\sn}{\operatorname{sn}}

\hypersetup{
  hidelinks,
  pdfauthor={},
  pdftitle={Treewidth and Separation number },
  pdfsubject={Treewidth and Separation number},
  pdfkeywords={treewidth, separation number, balanced separations}
}

\title{Treewidth and Separation Number:\\ An Optimal Linear Bound}
\author{Rong Chen \\
\\
Center for Discrete Mathematics,\ \ Fuzhou University\\
Fuzhou,\ \ P. R. China}

\begin{document}

\setlength{\abovedisplayskip}{6pt plus 2pt minus 2pt}
\setlength{\belowdisplayskip}{6pt plus 2pt minus 2pt}
\setlength{\abovedisplayshortskip}{0pt plus 2pt}
\setlength{\belowdisplayshortskip}{3pt plus 2pt minus 1pt}

\maketitle

\begin{abstract}

The \emph{separation number} $\sn(G)$ of a graph $G$ is the smallest nonnegative integer $a$ such that every subgraph $H$ of $G$ has a separation $(A,B)$ of {\em order} $|A\cap B|\le a$ for which both $A\setminus B$ and $B\setminus A$ contain at most $2|V(H)|/3$ vertices. Dvo\v{r}\'ak and Norin proved that $\tw(G)\le 15\sn(G)$ for every nonempty graph $G$. We improve this bound to $\tw(G)\le 4\sn(G)-2$. Moreover, we show that the coefficient 4 is asymptotically best possible: for every $k\ge 2$, there is a graph $G_k$ with $\tw(G_k)=4k$ and $\sn(G_k)=k+1$.


\end{abstract}

\footnotetext[1]{

Email: rongchen@fzu.edu.cn}

\section{Introduction}

All graphs in this paper are finite and simple. 
A {\em separation} $(A,B)$ of $G$ is a pair of subsets of $V(G)$ such that $A\cup B=V(G)$ and no edge of $G$ has one endpoint in $A\setminus B$ and the other in
$B\setminus A$. The {\em order} of $(A,B)$ is $|A\cap B|$.
Given a nonnegative weight function $w\colon V(G)\to\mathbb{R}_0^+$ and a vertex set $X\subseteq V(G)$, we define $w(X):=\sum_{v\in X}w(v)$. We say that a separation $(A,B)$ of a graph $G$ is
\emph{$w$-balanced} if $w(A\setminus B)\le\frac{2}{3}w(V(G))$ and
$w(B\setminus A)\le\frac{2}{3}w(V(G))$. A separation $(A,B)$ of a graph $G$ on
$n$ vertices is \emph{balanced} if $|A\setminus B|\le 2n/3$ and
$|B\setminus A|\le 2n/3$. Equivalently, it is $w$-balanced when $w$ ia the constant function equal to $1$. The {\em separation number} $\sn(G)$ is the smallest nonnegative integer $a$ such that every subgraph of $G$ has a balanced separation of order at
most $a$.

Let $\tw(G)$ denote the treewidth of a graph $G$. The relationship between the separation number and the treewidth has been studied since the work of Robertson and Seymour \cite{RS86}. They proved that if $\tw(G)\le t-1$, then every nonnegative weight function on $V(G)$ admits a $w$-balanced separation of $G$ of order at most $t$.
Conversely, they showed that if $G$ has a $w$-balanced separation of order
at most $a$ for every $\{0,1\}$-valued weight assignment $w$, then $\tw(G)<4a$.

Since $\tw(G')\le\tw(G)$ for every subgraph $G'$ of $G$, it follows that
$\sn(G)\le\tw(G)+1$ for every graph $G$. In this paper, we investigate whether a converse bound of a similar form holds. This converse is far from obvious, since the existence of small balanced separations in every subgraph does not directly imply the existence of a small $w$-balanced separation of the whole graph for every $\{0,1\}$-valued weight function $w$.
Dvo\v{r}\'ak and Norin~\cite{dvorak-norin} were the first to establish the following linear upper bound.

\begin{theorem}[Dvo\v{r}\'ak--Norin~\cite{dvorak-norin}] \label{thm:dvorak-norin}
Every nonempty graph $G$ satisfies $\tw(G)\leq15\sn(G)$.
\end{theorem}

Recently, Houdrouge, Miraftab and Morin~\cite{houdrouge-miraftab-morin} gave a constructive proof of the bound $\tw(G)\leq \frac{7915}{139}\sn(G)$. 
Prior to the work of Dvo\v{r}\'ak and Norin ~\cite{dvorak-norin}, Benjamini et al. \cite{BST12} gave a proof that there exists a function $f\colon\mathbb{N}\to\mathbb{N}$ such that $\tw(G)\le f(\sn(G))$ for
every graph $G$. Their proof uses the Excluded Grid Theorem of Robertson and Seymour \cite{RS86}. Using the current best version of the Excluded Grid Theorem \cite{CT21}, their proof shows that
$\tw(G)\le\sn^{56}(G)\cdot\log^{O(1)}\sn(G)$.


The proof in \cite{dvorak-norin} uses a haven derived from a bramble and a flow towards a suitable boundary set. Our proof of Theorem \ref{thm:upper-bound} follows the same general framework of haven and flow. In Section \ref{sec:upper-bound}, we improve the coefficient in Theorem \ref{thm:dvorak-norin} from 15 to 4, and show that this coefficient is asymptotically optimal in Section \ref{sec:lower-bound}. 


\begin{theorem}\label{thm:upper-bound}
Every nonempty graph $G$ satisfies $\tw(G)\leq4\sn(G)-2.$
\end{theorem}


\begin{theorem}\label{thm:lower-bound}
For every integer $k\ge2$, there exists a graph $G_k$ with \[  \tw(G_k)=4k  \ \text{and}\   \sn(G_k)=k+1.\]
\end{theorem}

Theorem \ref{thm:lower-bound} shows that the coefficient 4 in Theorem 1.2 is asymptotically best possible. It remains open whether the additive constant $-2$ can be improved. We conjecture that it cannot.



\section{Preliminaries}\label{sec:preliminaries}
The vertex set and edge set of a graph $G$ are denoted by $V(G)$ and $E(G)$, respectively. A \emph{separation} of $G$ is a pair $(A,B)$ of subsets of $V(G)$ such that $A \cup B = V(G)$ and there is no edge between $A \del B$ and $B \del A$.
The \emph{order} of $(A,B)$ is $|A \cap B|$.


A {\em tree decomposition} of a graph $G$ is a family $(B_x:x\in V(T))$ of
subsets of $V(G)$, called \emph{bags}, indexed by a tree $T$, such that every
vertex belongs to a bag, every edge has both endpoints in a bag, and for
each $v\in V(G)$ the nodes $x$ with $v\in B_x$ induce a connected subtree
of $T$. Its width is $\max_{x\in V(T)}|B_x|-1$, and $\tw(G)$ is the minimum
width of a tree decomposition of $G$.

A{\em bramble} in $G$ is a family $\mathcal B$ of nonempty connected subgraphs such that any two members intersect or have an edge between them.
A {\em hitting set} for $\mathcal B$ is a vertex set meeting every member of
$\mathcal B$. The {\em order} of $\mathcal B$ is the minimum size of a hitting set.

\begin{theorem}[Seymour--Thomas~\cite{seymour-thomas}]\label{thm:bramble-duality}
For every nonempty graph $G$, the maximum order of a bramble in $G$ is
$\tw(G)+1$.
\end{theorem}

Let $\mathcal B$ be a bramble of order $k\geq1$ in $G$. For each
$W\subseteq V(G)$ with $|W|<k$, some member of $\mathcal B$ avoids $W$ and all such members lie in a single component of $G\backslash W$, since they pairwise touch. Let $\beta(W)$ be the vertex set of this component. We call $\beta$ the {\em haven} defined by $\mathcal B$. Since $\beta(W)$ contains a member of $\mathcal B$, we have $|\beta(W)|\geq1$. Lemma \ref{lem:haven-properties} (1) and (2) are (\cite{dvorak-norin}, Observation 3).

\begin{lemma}\label{lem:haven-properties}
The haven $\beta$ defined by a bramble of order $k\geq1$ has the following properties.
\begin{enumerate}
\item[(1)] If $W\subseteq W'\subseteq V(G)$ and $|W'|<k$, then
  $\beta(W')\subseteq\beta(W)$.
\item[(2)] If $(C,D)$ is a separation of $G$ with $|C\cap D|<k$, and
  $W\subseteq V(G)$ satisfies $|W|<k$ and
  $\beta(W)\subseteq C\setminus D$, then
  $\beta(C\cap D)\subseteq C\setminus D$.
\item[(3)] For every $W\subseteq V(G)$ with $|W|<k$, we have
  $|W|+|\beta(W)|\geq k$.
\end{enumerate}
\end{lemma}

\begin{proof}
For (3), since $W\cup\beta(W)$ is a hitting set for $\mathcal B$ and $\beta$ is defined by a bramble of order $k\geq1$, we have $|W|+|\beta(W)|\geq k$.
\end{proof}

Let $\vec{G}$ be a directed graph and let $g\colon V(\vec{G})\to \mathbb{R}_0^+$
be a function specifying the amount of flow entering its vertices. The value
$g(V(\vec{G})) := \sum_{v\in V(\vec{G})} g(v)$ is called the \emph{total supply}.
Let $W\subseteq V(\vec{G})$ be a set of sinks in $\vec{G}$, where a \emph{sink} is a vertex with no outgoing edges.
A \emph{flow towards $W$} in $\vec{G}$ with supply $g$ is a function
$f\colon E(\vec{G})\to\mathbb{R}_0^+$ satisfying the flow conservation equation
\[ g(v) + \sum_{e=uv\in E(\vec{G})} f(e)= \sum_{e=vz\in E(\vec{G})} f(e)
   \qquad\text{for every }v\in V(\vec{G})\setminus W.\]
For a fixed flow $f$ and a vertex $v\in V(\vec{G})$, define $\operatorname{in}(v)$ to be the sum of the flow values on the edges into $v$. The \emph{congestion at $v$} is defined as
$\operatorname{in}(v)+g(v)$. The \emph{congestion of a flow $f$} is the maximum
congestion over all the vertices of the graph. Summing the conservation equations shows that the total supply $\sum_{v\in V(\vec{G})} g(v)$ equals the sum of congestions of vertices of $W$. Let the \emph{support} of $f$, denoted by $\operatorname{supp}(f)$, consist of $\vec{G}[W]$ together with the set of edges of $\vec{G}$ carrying a non-zero flow value. By the conservation equations, each vertex with a positive supply is in $\operatorname{supp}(f)$. 

Let $G$ be an undirected graph with $W\subseteq V(G)$. The graph
$\vec{G}_W$ is obtained from $G$ by replacing each edge with a pair of
oppositely directed edges and deleting the edges directed away from the
vertices of $W$. That is, no arc links two vertices of $W$. A \emph{$W$-cloud} in $G$ is a flow in $\vec{G}_W$ towards $W$ such that the supply at each vertex is at most $1$. For a positive integer $m$, let $\nu_m(G,W)$ be the maximum total supply
of a $W$-cloud in $G$ with congestion at most $m$.

Lemma \ref{lem:integral-support} (1) is (\cite{dvorak-norin}, Observation 5) without proof. For completeness, we give a proof here.

\begin{lemma}\label{lem:integral-support}
Let $G$ be a graph with $W\subseteq V(G)$, and let $m\geq1$ be an integer. Then
\begin{itemize}
\item[(1)] $$\nu_m(G,W)= \min_{\substack{(A,B)\text{ a separation of }G\\W\subseteq B}}
  \bigl(|B\setminus A|+m|A\cap B|\bigr).$$
\item[(2)] There is a subgraph $H$ of $G$ with $W\subseteq V(H)$, such that
\begin{itemize}
\item[(2.1)] $|V(H)|=\nu_m(G,W)$;
\item[(2.2)] when $W=\emptyset$, we have $V(H)=\emptyset$; and
\item[(2.3)] for any separation $(A,B)$ of $G$, we have \[|A\cap V(H)|\leq m\,|A\cap(W\cup B)\cap V(H)|.\]
\end{itemize}
\end{itemize}
\end{lemma}

\begin{proof}
First, we prove (1). Split every vertex $v\in V(\vec{G}_W)$ into an input $v_{\rm in}$ and an output
$v_{\rm out}$, adding an arc from $v_{\rm in}$ to $v_{\rm out}$ 
of capacity $m$. Add a source $s$ with an arc of capacity $1$ to each input.
For every arc $uv\in E(\vec{G}_W)$, add an arc from $u_{\rm out}$ to $v_{\rm in}$ of infinite capacity. Join $w_{\rm out}$ to a new sink $z$ by an infinite-capacity arc for every $w\in W$.
Let $N$ denote this network. Note that
the flow in $N$ on $sv_{\rm in}$ is the supply at $v$, and the flow on $v_{\rm in}v_{\rm out}$ is its congestion. Since $W$-clouds in $G$ with congestion at most $m$ are in one to one correspondence with the flows in $N$, the maximum flow value in $N$ is $\nu_m(G,W)$. Hence, by the max-flow min-cut theorem, to prove (1), it suffices to show that the capacity of a minimum $s$--$z$ cut in $N$ is $$\min_{\substack{(A,B)\text{ a separation of }G\\W\subseteq B}}  \bigl(|B\setminus A|+m|A\cap B|\bigr).$$ 

Let $(S,S^c)$ be a minimum $s$--$z$ cut in $N$ with $s\in S$ and $z\in S^c$. Evidently,
$(S,S^c)$ contains no infinite-capacity arc, so

\ \ (a) $w_{\rm out}\in S^c$ for each $w\in W$

\noindent since the arc $w_{\rm out}z$ has infinite capacity. Assume $v_{\rm in}\in S^c$ for some $v\in V(G)$. If  $v_{\rm out}\in S$, then ${\rm Cap}(S\cup \{v_{\rm in}\},S^c-\{v_{\rm in}\})={\rm Cap}(S,S^c)-1$,  contradicting minimality, where ${\rm Cap}(S,S^c)$ denotes the capacity of the cut $(S,S^c)$. So  $v_{\rm out}\in S^c$.
That is, for some $v\in V(G)$,

\ \ (b) if $v_{\rm in}$ and $v_{\rm out}$ are not in the same side of $(S,S^c)$, then $v_{\rm in}\in S$ and $v_{\rm out}\in S^c$.

Let $A$ consist of vertices in $G$ whose input is in $S$, and let $B$ consist of vertices in $G$ whose output is in $S^c$, implying $W\subseteq B$ by (a). That is, $A\cap B$ consists of vertices whose inputs are in $S$ and whose outputs are in $S^c$, and  $A\setminus B$ consists of vertices in $G$ whose input and output are in $S$, and $B\setminus A$ consists of vertices in $G$ whose input and output are in $S^c$. Then $V(G)=A\cup B$ by (b). Moreover, since $(S,S^c)$ contains no infinite-capacity arc and an edge between $A\setminus B$ and $B\setminus A$ would give an infinite-capacity arc across $(S,S^c)$, we have that $(A,B)$ is a separation of $G$ with $W\subseteq B$ by (a). 
Hence, ${\rm Cap}(S,S^c)=|B\setminus A|+m|A\cap B|$ by the definition of $N$. Conversely, every separation occurring in (1) defines an $s$--$z$ cut in $N$ with this capacity. This proves (1).


Secondly, we prove (2). Since all capacities in $N$ are integers, $N$ has an integral maximum flow $f'$ of value $\nu_m(G,W)$. Let $f$ be the $W$-cloud in $\vec{G}_W$ corresponding to $f'$. Evidently, $f$ is an integral $W$-cloud whose total supply is $\nu_m(G,W)$. Among all integral maximum $W$-clouds in $\vec{G}_W$, we may assume that $f$ is chosen minimizing the sum of its values on graph arcs. If $\operatorname{supp}(f)$ has a directed cycle $C$, subtracting one unit around $C$ would preserve the total supply and all constraints while decreasing this sum. Hence, $\operatorname{supp}(f)$ has no directed cycle. In particular, the underlying graph of $\operatorname{supp}(f)$, denoted by $H$, is a subgraph of $G$ and $V(H)=V(\operatorname{supp}(f))$. Then $W\subseteq V(H)$ by the definition of supply.

Suppose that a vertex $v\in V(H)$ has an in-neighbour in $\operatorname{supp}(f)$ and satisfies $g(v)=0$. Since $\operatorname{supp}(f)$ has no directed cycle, there is a directed path in $\operatorname{supp}(f)$ from a vertex $x$ to $v$ with $g(x)=1$. Subtract one unit along this path, decrease the supply at $x$ by one, and increase the supply at $v$ by one. The total supply is unchanged. Congestion decreases at $x$ and at the internal vertices of the path, while congestion at $v$ is unchanged.
The resulting integral flow is feasible and has a smaller sum of graph arc values, a contradiction. So $g(v)=1$ as $f$ is an integral flow.

We claim that $g(v)=1$ for each vertex $v\in V(H)$. Assume to the contrary that $g(v)=0$ for some $v\in V(H)$. By the result proved in the last paragraph, $v$ has no in-neighbour in $\operatorname{supp}(f)$, implying that the congestion at $v$ is 0, so no arc incident with $v$ carries a non-zero flow value. Hence, $v\in W$ by the definition of support. 
Adding one unit of supply at $v$ would then respect the congestion bound $m\geq1$ and increase the total supply, a contradiction.

Since each vertex in $V(H)-W$ with a positive supply is in $\operatorname{supp}(f)$ by the conservation equations, $g(v)=0$ for each vertex $v\notin V(H)$.
Moreover, by the claim proved in the last paragraph, the total supply is $|V(H)|$, implying  $|V(H)|=\nu_m(G,W)$, proving (2.1). When $W=\emptyset$, summing conservation equations gives total supply $0$, and the empty support has all the required properties, proving (2.2).

Decompose $f$ into unit paths from supply vertices to $W$, allowing a path of length zero for supply at a terminal. Every path starting in $A\cap V(H)$ meets $A\cap(W\cup B)\cap V(H)$. Indeed, it either ends at a terminal in $A\cap W$, or passes through $A\cap B\cap V(H)$ before reaching a terminal outside $A$. Charge its unit of supply to the first vertex of $A\cap(W\cup B)\cap V(H)$ on the path. At any vertex the number of charged paths is at most its congestion, hence at most $m$. Since $g(v)=1$ for each vertex $v\in V(H)$, there is one supply unit starting at every vertex of $A\cap V(H)$, which proves (2.3).
\end{proof}

\section{Proof of Theorem \ref{thm:upper-bound}}\label{sec:upper-bound}
Let $G$ be a non-empty graph. In this section, set \[a:=\sn(G),\ \ \ k:=\tw(G)+1.\]
Since $\sn(G)\leq \tw(G)+1$, we have that $k\ge a$ by Theorem \ref{thm:bramble-duality}.
Let $\beta$ be the haven defined by a bramble of order $k$ in $G$. For each integer $j\ge1$, define
\begin{eqnarray*}
\Psi_j(W)&:=&|\beta(W)|+j|W|\quad\text{for any}\ W\subseteq V(G)\text{ with }|W|<k,\\
  p_j&:=&\min_{|W|<k}\Psi_j(W).
\end{eqnarray*}
Note that $\Psi_j(\emptyset)=|\beta(\emptyset)|$ for all integers $j\geq1$.
Since $\Psi_{j+1}(W)=\Psi_j(W)+|W|$, we have $p_{j+1}\geq p_j$. Moreover, since $p_1\geq k$ by Lemma~\ref{lem:haven-properties} (3), we have
\begin{equation}\label{eq:P1}
  k\leq p_1\leq p_2\leq p_3\leq\ldots\leq|\beta(\emptyset)|.
\end{equation}
For each integer $r\geq0$ with $ra<k$, set
\[ \psi_j(r):=\min_{|W|+ra<k}\Psi_j(W).\]
Since its domain is finite and contains the empty set, all $\psi_j(r)$ exist. Evidently,
\begin{equation}\label{eq:phi}
  p_j=\psi_j(0)\leq\psi_j(1)\leq\psi_j(2)\leq\psi_j(3)\leq\ldots\leq|\beta(\emptyset)|.
\end{equation}
The condition $|W|+ra<k$ leaves room to add at most $a$ vertices at each of $r$ subsequent steps.

\begin{lemma}\label{minimizer W}
When $a<k$, there is an integer $m$ with $1\leq m\le|\beta(\emptyset)|$ such that some global
minimizer $W$ of $\Psi_m$ satisfies $|W|+a< k$.
\end{lemma}
\begin{proof}
For any non-empty subset $W$ of $V(G)$, since
$$\Psi_{|\beta(\emptyset)|}(W)\ge|\beta(\emptyset)|\cdot|W|+1>|\beta(\emptyset)|
=\Psi_{|\beta(\emptyset)|}(\emptyset),$$ the set $\emptyset$ is a
global minimizer of $\Psi_{|\beta(\emptyset)|}$, so the lemma holds as $|\emptyset|+a=a<k$.
\end{proof}

\begin{lemma}\label{minimizer m}
Let $m\geq1$ be an integer. If $k\geq 4a$ and every global minimizer $W$ of $\Psi_m$ satisfies $|W|+a\ge k$, then $p_m-3ma\geq a$.
\end{lemma}
\begin{proof}
By Lemma~\ref{lem:haven-properties} (3) and the fact that $k\geq 4a$, we have
\[ p_{m}=|\beta(W)|+m|W|\geq k+(m-1)|W|\geq mk-(m-1)a\geq3ma+a,\]
so $p_m-3ma\geq a$.
\end{proof}

\begin{lemma}\label{lem:budget-recurrence}
For all integers $j\geq1$ and $r\geq1$ with $ra<k$, we have
\begin{equation}\label{eq:budget-recurrence}
  \psi_j(r)\geq\frac32\psi_j(r-1)-\frac32ja.
\end{equation}
Consequently, for every integer $r\geq0$ with $ra<k$,
\begin{equation}\label{eq:iterated-recurrence}
  \psi_j(r)\geq3ja+\left(\frac32\right)^r(p_j-3ja).
\end{equation}
\end{lemma}

\begin{proof}
Choose $W\subseteq V(G)$ such that $|W|+ra<k$ and $\psi_j(r)=\Psi_j(W)$. Since its domain is finite and contains the empty set, such $W$ exists.  Set $ R:=\beta(W)$. That is, $\psi_j(r)=\Psi_j(W)=|R|+j|W|$.

\begin{claim}\label{clm:boundary-flow}
Every vertex of $W$ has a neighbor in $R$ and 
$  \nu_j(G\backslash R,W)=j|W|.$
\end{claim}
\begin{proof}
Suppose to the contrary that some $w\in W$ has no neighbor in $R$. Since $G[R]$ is a component of $G\backslash W$, the subgraph $G[R]$ is also a component of $G\backslash (W\setminus\{w\})$. Moreover, since $R=\beta(W)$, we have $\beta(W\setminus\{w\})=R$ by Lemma \ref{lem:haven-properties} (1), so $\Psi_j(W\setminus\{w\})=\Psi_j(W)-j<\Psi_j(W)$ as $j\geq1$, contradicting the choice of $W$. Hence, each vertex of $W$ has a neighbor in $R$.

Now we prove that $\nu_j(G\backslash R,W)=j|W|.$ By Lemma~\ref{lem:integral-support} (1), it suffices to show
\[ j|W|=\min_{\substack{(A,B)\text{ a separation of }G\backslash R\\W\subseteq A}}\bigl(|A\setminus B|+j|A\cap B|\bigr).\]
Since $j|W|=|A\setminus B|+j|A\cap B|$ when $(A,B)=(W,V(G\backslash R))$, it suffices to show that

\ \   (a) $j|W|\leq |A\setminus B|+j|A\cap B|$ for each separation $(A,B)$ of $G\backslash R$ with $W\subseteq A$.

\noindent When $|A\cap B|\geq |W|$, (a) obviously holds. So we may assume that $|A\cap B|< |W|$. Since $R$ is a component of $G\del W$ and $W\subseteq A$, we have that $(A\cup R,B)$ is a separation of $G$ with order $|A\cap B|$. By Lemma~\ref{lem:haven-properties} (2), $\beta(A\cap B)\subseteq R\cup(A\setminus B).$
The minimizing property of $W$ now implies
\[  |R|+j|W|
  \leq|\beta(A\cap B)|+j|A\cap B|
  \leq |R|+|A\setminus B|+j|A\cap B|,\]
so (a) holds.
\end{proof}

By Claim~\ref{clm:boundary-flow} and Lemma~\ref{lem:integral-support} (2) (with $G$ replaced by $G\backslash R$), there is a subgraph $H$ of $G\backslash R$ with $W\subseteq V(H)$ and $|V(H)|=j|W|$ and satisfying Lemma~\ref{lem:integral-support} (2.2) and (2.3). Set $F:=H\cup G[W\cup R].$ Since $(V(H), R)$ is a partition $V(F)$, we have
\begin{equation}\label{eq:auxiliary-order}
  |V(F)|=|V(H)|+|R|=j|W|+|R|=\psi_j(r).
\end{equation}
When $W=\emptyset$, since $V(H)=\emptyset$ by Lemma~\ref{lem:integral-support} (2.2), we have $F=G[R]$. Since $F$ is a subgraph of $G$, there is a balanced separation $(A,B)$
of $F$ with $|A\cap B|\leq a$. 

\begin{claim}\label{clm:small-oriented-side}
After interchanging $A$ and $B$ if necessary, we have $|A|\leq ja+\psi_j(r)-\psi_j(r-1)$.
\end{claim}
\begin{proof}
Suppose first that $R\subseteq B\setminus A$. Since $N_G(R)\subseteq W$, we have $W\subseteq B$, so $A\cap(W\cup B)\cap V(H)\subseteq A\cap B$. Moreover, by Lemma~\ref{lem:integral-support} (2.3) and (\ref{eq:phi}), 
\[  |A|=|A\cap V(H)|\leq j|A\cap B|\leq ja\leq ja+\psi_j(r)-\psi_j(r-1),\]
so Claim \ref{clm:small-oriented-side}  holds. Hence, by symmetry, we may assume that $R\nsubseteq B\setminus A$ and $R\nsubseteq A\setminus B$.

Set $T:=A\cap B\cap R$ and $W':=W\cup T$. Since $r\geq1$, we have $ |W'|\leq |W|+a\leq |W|+ra<k.$
Set $R':=\beta(W')$. By Lemma~\ref{lem:haven-properties} (1), $R'$ is contained in $R\setminus T$ and is disjoint from $A\cap B$. Since $F[R']$ is a connected subgraph of $F\backslash (A\cap B)$,  by interchange $A$ and $B$ if necessary, we may assume that $R'\subseteq B\setminus A$. Set
\[  W'':=(W\cap B)\cup T, \qquad  \delta:=|T|-|W\setminus B|.\]
Then $|W''|=|W|+\delta$ and $\delta\leq|T|\leq a$, so
\begin{equation}\label{eq:restored-budget}
  |W''|+(r-1)a\leq |W|+ra<k.
\end{equation}

We claim that $R'= \beta(W'')$. 
Since $F[R']$ is a connected subgraph of $F\backslash (A\cap B)$ and $R'\subseteq (B\setminus A)\cap R$, no edge in $F$ links $R'$ and $W\del B$. Moreover, since $W''\subseteq W'$ and $R'=\beta(W')$, the graph $G[R']$ is a component of $G\del W''$ because deleting $W\del B$ removes only vertices having no neighbour in $R'$, so $R'= \beta(W'')$. This proves the claim.

Since $R'= \beta(W'')$ and $|W''|=|W|+\delta$, by (\ref{eq:restored-budget}) we have $$\psi_j(r-1)\leq \beta(W'')|+j|W''|=|R'|+j(|W|+\delta).$$ Combining this with the fact that $\psi_j(r)=j|W|+|R|$, we have
\[|R|-|R'|\leq \psi_j(r)-j|W|+j(|W|+\delta)-\psi_j(r-1)=j\delta+\psi_j(r)-\psi_j(r-1).\]
Moreover, since  $R'\subseteq (B\setminus A)\cap R$, by Lemma \ref{lem:integral-support} (2.3),
\begin{eqnarray*}
|A|\leq |R|-|R'|+|V(H)\backslash B| &\leq& |R|-|R'|+j(|(A\cap B)\backslash T|+|W\backslash B|)\\
  &\leq& |R|-|R'|+j(a-|T|+|W\backslash B|)\\
  &=& |R|-|R'|+j(a-\delta)\\
  &\leq& ja+\psi_j(r)-\psi_j(r-1).
\end{eqnarray*}
This proves Claim \ref{clm:small-oriented-side}.
\end{proof}

Since $(A,B)$ is balanced in $F$ and $|V(F)|=\psi_j(r)$ by \eqref{eq:auxiliary-order}, we have $|A|\geq \frac13\psi_j(r)$.
Combined with Claim~\ref{clm:small-oriented-side}, we have $\psi_j(r)\leq 3ja+3\bigl(\psi_j(r)-\psi_j(r-1)\bigr),$ so \eqref{eq:budget-recurrence} holds and
 \[  \psi_j(r)-3ja\geq\frac32\bigl(\psi_j(r-1)-3ja\bigr).\]
Since $\psi_j(0)=p_j$, we have that \eqref{eq:iterated-recurrence} holds.
\end{proof}

\begin{corollary}\label{cor:interior-minimizer}
Let $m\geq1$ be an integer. If some global minimizer $W$ of $\Psi_m$
satisfies $|W|+a<k$, then $p_m\leq3ma$.
\end{corollary}
\begin{proof}
The set $W$ is admissible for $\psi_m(1)$, so
$\psi_m(1)=\psi_m(0)=p_m$. Applying
\eqref{eq:budget-recurrence} with $r=1$ gives $ p_m\geq\frac32p_m-\frac32ma,$
so $p_m\leq3ma$.
\end{proof}

We can now prove Theorem \ref{thm:upper-bound}. For convenience, we restate it here.
\setcounter{theorem}{1}
\renewcommand{\thetheorem}{1.\arabic{theorem}}
\begin{theorem}
Every nonempty graph $G$ satisfies $\tw(G)+1\leq4\sn(G)-1.$
\end{theorem}
\begin{proof}
Set $a:=\sn(G)$ and $k:=\tw(G)+1$. Since a one-vertex subgraph requires a balanced separation of order one, $a\geq1$. Suppose for a contradiction that $k>4a-1$. Since $k$ and $a$ are integers, $k\geq4a$.
By Theorem~\ref{thm:bramble-duality}, $G$ has a bramble of order $k$. Let $\beta$ be the haven defined by this bramble. 

By Lemma \ref{minimizer W}, there is an integer $1\leq m\le|\beta(\emptyset)|$ such that some global
minimizer $W$ of $\Psi_m$ satisfies $|W|+a< k$. Then $ p_m\leq3ma$ by Corollary~\ref{cor:interior-minimizer}, implying
\begin{equation}\label{eq:W-potential-upper}
  \Psi_{m-1}(W)=\Psi_m(W)-|W|=p_m-|W|\leq3ma-|W|,
\end{equation}
Moreover, since $4a\leq k\leq p_1\leq p_m\leq3ma$ by \eqref{eq:P1}, we have $m\geq2$. Among all $m$ satisfying the conclusion  of Lemma \ref{minimizer W}, choose the smallest.
By the minimality of $m$, every global minimizer $Z$ of $\Psi_{m-1}$ satisfies $|Z|\geq k-a$. Then it follows from Lemma~\ref{minimizer m} that
\begin{equation}\label{eq:Pj-lower}
 p_{m-1}-3(m-1)a\ge a.
\end{equation}

Since $|W|+a<k$, the set $W$ is admissible for $\psi_{m-1}(1)$.
By \eqref{eq:iterated-recurrence}, \eqref{eq:W-potential-upper} and \eqref{eq:Pj-lower},
\[3ma-|W|\ge\Psi_{m-1}(W)\geq \psi_{m-1}(1)  \geq3(m-1)a+\frac32a  =\left(3m-\frac32\right)a,\]
so $|W|\leq3a/2$. Since $|W|+2a\leq\frac72a<4a\leq k,$
the set $W$ is also admissible for $\psi_{m-1}(2)$. By \eqref{eq:iterated-recurrence}, \eqref{eq:W-potential-upper} and \eqref{eq:Pj-lower} again,
\[3ma-|W|\ge\Psi_{m-1}(W)\geq \psi_{m-1}(2)  \geq3(m-1)a+\frac94a  =\left(3m-\frac34\right)a,\]
implying $|W|\leq3a/4$, so $ |W|+3a\leq\frac{15}{4}a<4a\leq k$.
Since $W$ is admissible for $\psi_{m-1}(3)$, by \eqref{eq:iterated-recurrence}, \eqref{eq:W-potential-upper} and \eqref{eq:Pj-lower} again
\[ 3ma-|W|\ge\Psi_{m-1}(W)\geq\psi_{m-1}(3)  \geq3(m-1)a+\frac{27}{8}a  =\left(3m+\frac38\right)a>3ma,\]
a contradiction. 
\end{proof}

\section{Proof of Theorem~\ref{thm:lower-bound}}\label{sec:lower-bound}
Let $G$ be a graph with $n$ vertices. For any $S\subseteq V(G)$, if each component of $G\del S$ has at most ${2n}/3$ vertices, we say that $S$ is a \emph{balanced cut} of $G$.
\setcounter{theorem}{0}
\renewcommand{\thetheorem}{4.\arabic{theorem}}
\begin{lemma}\label{lem:component-separator}
Let $G$ be a graph with $n\geq1$ vertices and let $S\subseteq V(G)$. 
\begin{enumerate}
\item[(1)] $S$ is a balanced cut of $G$ if and only if there is a balanced separation $(A,B)$  of $G$ with $A\cap B=S$.
\item[(2)] If all components of $G\del S$ with size larger than ${2n}/3$ are trees, then $G$ has a balanced cut of size at most $|S|+1$.
\end{enumerate}
\end{lemma}
\begin{proof}
First we consider (1). The ``if" part is obvious. So it suffices to show that the ``only if" part is true. Assume that $S$ is a balanced cut of $G$. It suffices to partition the components of $G\del S$ into two collections whose unions each have at most $2n/3$ vertices. If some component has between $n/3$ and $2n/3$ vertices, put it on one side and all other components on the other. Otherwise, every component has fewer than $n/3$ vertices. If their total size is less than $n/3$, put all of them on one side. In the remaining case, accumulate components until their union first has at least $n/3$ vertices. This union has fewer than $2n/3$ vertices, and its complement in $V(G)\setminus S$ has at most $2n/3$ vertices. Adding $S$ to both sides gives the desired separation.

Note that there is at most one component of $G\del S$ with more than $2n/3$ vertices. If there is none,
then $S$ is a balanced cut of $G$ by definition. 
So we may assume that there is a unique component $T$ of $G\del S$ with $|V(T)|>2n/3$. Since $T$ is a tree, there is a vertex $v\in V(T)$ such that each component of $T\del v$ has at most $|V(T)|/2$ vertices. Since $|V(T)|/2\leq n/2$, each component of $G\del (S\cup\{v\})$ has at most $2n/3$ vertices, so $S\cup\{v\}$ is a balanced cut of $G$.
\end{proof}

The \emph{strong separation number} of a graph $G$, denoted by $\operatorname{sn}^{*}(G)$, is the smallest integer $a$ such that every connected subgraph of $G$ has a balanced cut of size at most $a$. Then $\operatorname{sn}^{*}(G)\le\operatorname{sn}(G)$ by Lemma \ref{lem:component-separator} (1).
\begin{lemma}\label{lem:sn=sn*}
$\operatorname{sn}(G)=\operatorname{sn}^{*}(G)$.
\end{lemma}
\begin{proof}
Since $\operatorname{sn}^{*}(G)\le\operatorname{sn}(G)$, to prove the lemma, it suffices to show that  each disconnected subgraph $H$ of $G$ has a balanced cut of size at most $\operatorname{sn}^{*}(G)$. When each component of $H$ has size at most $\frac{2}{3}|V(H)|$, the empty set is a balanced cut of $H$. So we may assume that some component $C$ of $H$ with $|V(C)|>\frac{2}{3}|V(H)|$. Since $C$ is a connected subgraph of $G$, there is a set $X\subseteq V(C)$ with $|X|\le\operatorname{sn}^{*}(G)$ such that each component of $C\setminus X$ has size at most $\frac{2}{3}|V(C)|$. Since except $C$, all other components of $H$ have size less than $\frac{1}{3}|V(H)|$, each component of $H\setminus X$ has size at most $\frac{2}{3}|V(H)|$.
\end{proof}

Theorem~\ref{thm:lower-bound} is an immediate consequence of Theorem \ref{Ht} by setting $G_k:=H_{4k+1}$.
\begin{theorem}\label{Ht}
For every integer $t\geq8$ with $t\equiv0,1\pmod4$, there exists a subdivision $H_t$ of $K_t$ such that \[  \tw(H_t)=t-1  \ \text{and}\   \sn(H_t)=\lfloor \frac t4\rfloor+1.\]
\end{theorem}
\begin{proof}
Fix an integer $t\geq8$, and label the vertices of $K_t$ by $\{1, 2, \dots, t\}$. Set $\ell_r:=(10t^3)^r$ for each $0\leq r\leq t-1.$
For each edge $ij$ of $K_t$ with $i<j$, replace $ij$ by a path $P_{ij}$ of length $\ell_i$, using new internal vertices for each path. Let $H_t$ be the resulting subdivision. The original vertices of $K_t$ are called {\em branch vertices} in $H_t$, and the {\em rank}  of $P_{ij}$ is $i$. Since $\tw(K_t)=t-1$, we have $\tw(H_t)=t-1$. 

Since $t/4+1/4\leq \sn(H_t)$ by Theorem \ref{thm:upper-bound}, we have $\sn(H_t)=\lfloor \frac t4\rfloor+1$ if $\sn(H_t)\le\lfloor \frac t4\rfloor+1$. Hence, to prove Theorem \ref{Ht}, it suffices to show that $\sn(H_t)\le\lfloor \frac t4\rfloor+1$.
Moreover, to prove $\operatorname{sn}(H_t)\leq\lfloor\frac t4\rfloor+1$, by Lemmas \ref{lem:component-separator} (1) and \ref{lem:sn=sn*}, it suffices to show that any non-empty connected subgraph $H$ of $H_t$ has a balanced cut of size at most $\lfloor\frac t4\rfloor+1$. 

Call an original path $P_{ij}$ {\em complete} in $H$ if  $P_{ij}\subseteq H$. Evidently, when no original path is complete, $H$ is a forest, so $H$ has balanced cut of size one by Lemma \ref{lem:component-separator} (2) (in this case, choose $S=\emptyset$).
We may therefore assume that some original path is complete in $H$. Let $i$ be the maximum rank of a complete original path in $H$. Then $1\leq i\leq t-1$ and the branch vertex $i$ belongs to $H$. Set $n:=|V(H)|$. Evidently, $ n\geq \ell_i$. Since $H\del ([i-1]\cap V(H))$ 
is a forest, $H$ has a balanced cut of size at most $i$ by Lemma \ref{lem:component-separator} (2) again. Hence, when $i\le \lfloor\frac t4\rfloor+1$, we are done. So we may assume that $i>\lfloor\frac t4\rfloor+1\ge3$ as $t\ge8$.

Call the branch vertices $1,\ldots,i-1$ {\em early} and the branch vertices $i+1,\ldots,t$ {\em late}.
Let $X$ consist of the early branch vertices in $H$ and all vertices of $H$ internal to original paths of rank less than $i$. There are at most $\binom{t}{2}$ original paths of rank less than $i$,
each of length at most $\ell_{i-1}$. Moreover, since $ n\geq \ell_i=10t^3\ell_{i-1}$,
\begin{equation}\label{eq:short-part}
  |X\cup\{i\}|=|X|+1  \leq t-1+\binom{t}{2}\ell_{i-1}  < t^2\ell_{i-1} =\frac{\ell_i}{10t}  \leq\frac{n}{10t}.
\end{equation}

Let $C$ be a component of $H\setminus(X\cup\{i\})$. Since 
$i$ is the maximal rank of a complete original path in $H$, the graph 
$C$ contains at most one late branch vertex, and $C$ is a path consisting of internal vertices of an original path when it does not contain a late branch vertex. That is, no matter which case happens, $C$ is a tree. 
When $C$ contains a late branch vertex, let $v_C$ denote the branch vertex in $C$; otherwise let $v_C$ denote a vertex in $C$ with the sizes of components of $C\del v_C$ as equal as possible. We say that $v_C$ is the \emph {bundle} of $C$. Since every edge between $X$ and $C$ is incident with $v_C$, 

\ \ \ (a) there is no edge between $X$ and $C\del v_C$.

When there is a component $C$ of $H\setminus(X\cup\{i\})$ with $|V(C)|\ge\frac {2n}3$, since $C$ is a tree, there is a balanced cut of $H$ with size at most 3 containing $\{i,v_C\}$ by (a) and Lemma \ref{lem:component-separator} (2) (in this case, choose $S=\{i,v_C\}$). Since $\lfloor\frac t4\rfloor+1\ge3$, under this case, we are done. Hence, we may assume that each component of $H\setminus(X\cup\{i\})$ has size less than ${2n}/3$.


Set $k:=\lfloor \frac{t-i}{3}\rfloor+1$ and $\lambda:= \frac k{t-i}$. Since $t-i$ is a positive integer, $k\ge\frac{t-i}{3}+\frac13$, so $\lambda\ge\frac13+\frac{1}{3(t-i)}$. Hence, by
\eqref{eq:short-part},
\begin{equation}\label{eq:big-part}
  \lambda(n-|X|-1)\ >\ \lambda\left(n-\frac{n}{10t}\right)
  \ \ge\ \frac{n}{3}\left(1+\frac{1}{t-i}\right)\left(1-\frac{1}{10t}\right)>\frac{n}{3}.
\end{equation}

Let $\mathcal C_1$ be the collection of components of $H\setminus (X\cup\{i\})$ containing a late branch vertex, and let $\mathcal C_0$ be the collection of the remaining components. Since each component in $\mathcal C_1$ contains exactly one late branch vertex, we have $|\mathcal C_1|\le t-i$. Moreover, no edge of $H$ joins $X$ to a component in $\mathcal C_0$. For any $j\in\{0,1\}$, set $w_j:=\sum_{C\in\mathcal C_j}|V(C)|. $ Then $w_0+w_1=n-|X|-1$.
Choose the $\min\{k,|\mathcal C_1|\}$ largest components in $\mathcal C_1$, and let $\mathcal D$ be the collection of these components. Set $S:=\{v_C:C\in\mathcal D\}$ and $w:=\sum_{C\in\mathcal D}|V(C)|.$
Since $|\mathcal C_1|\le t-i$, we have $w\ge\lambda w_1$, so
\[w_0+w\ge\lambda(w_0+w_1)=\lambda(n-|X|-1)>\frac n3. \eqno{(11)}\]

Since $i\ge\lfloor \frac t4\rfloor+2$ and $t\equiv0,1\pmod4$, we have 
$\lfloor \frac{t-i}{3}\rfloor\le \lfloor\frac t4\rfloor-1$, so $|S\cup\{i\}|\le k+1\leq\lfloor \frac t4\rfloor+1$. Hence, 
it suffices to show that $S\cup\{i\}$ is a balanced cut of $H$. We claim that $S\cup\{i\}$ is a balanced cut of $H$. Let $C$ be a component of $H\setminus(S\cup\{i\})$. If $V(C)\cap X=\emptyset$, then $C$ is contained in a component of $H\setminus(X\cup{i})$, so $|V(C)|<2n/3$. So we may assume that $V(C)\cap X\neq\emptyset$. Since $H$ is connected, $N_H(C')=\{i\}$ for any component $C' \in\mathcal C_0$,
so the component $C$ contains no vertex of any component in $\mathcal C_0\cup\mathcal D$ by (a). Thus, $|V(C)|\le n-1-w_0-w<\frac{2n}{3}$ by (11). This proves the claim.
\end{proof}



\section{Acknowledgments}
This research was partially supported by grants from the National Natural Sciences Foundation of China (No. 12671400) and grants from the Natural Sciences Foundation of Fujian Province (No. 2025J01486).

	\vspace{5pt}
\noindent$\textbf{Disclosure of AI use}$

We used AI in two ways. First, based on the approach of \cite{dvorak-norin}, ChatGPT 6pro suggested an improvement of the upper bound in Theorem \ref{thm:dvorak-norin} to $\tw(G)\le 4\sn(G)-2$, together with examples showing that $\tw(G)\ge(4-\varepsilon)\sn(G)$ is attainable for every $\varepsilon>0$. 
We subsequently developed a stronger construction and proved Theorem \ref{Ht}, which implies Theorem \ref{thm:lower-bound} and establishes the asymptotic optimality of the coefficient 4.

Second, we used DeepSeek Harness and ChatGPT pro to proofread the English and to check the logic of the argument. We verified every mathematical claim
independently, rewrote and reorganized all AI-assisted material in our own
notation, and take full responsibility for the content of the paper.


\end{document}